\documentclass[11pt]{article}
\usepackage{amsmath,amsthm,amssymb,amsfonts}
\usepackage{amsfonts,amssymb,color}

\def\G{\overline{G}}

\newtheorem{theorem}{Theorem}[section]
\newtheorem{proposition}[theorem]{Proposition}

\newtheorem{corollary}[theorem]{Corollary}
\newtheorem{remark}{Remark}
\newtheorem{conjecture}{Conjecture}[section]
\numberwithin{equation}{section}

\title{Graph toughness from Laplacian eigenvalues}

\author{Xiaofeng Gu\thanks{
University of West Georgia, Carrollton, GA 30118, USA. Email: {\tt xgu@westga.edu}}
\ \ and\ \ Willem H. Haemers\thanks{Tilburg University, Tilburg, The Netherlands. Email: {\tt haemers@uvt.nl}}
}

\begin{document}
\date{}
\maketitle
\noindent

\begin{abstract}
The toughness $t(G)$ of a graph $G=(V,E)$ is defined as $t(G)=\min\{\frac{|S|}{c(G-S)}\}$,
in which the minimum is taken over all $S\subset V$ such that $G-S$ is disconnected,
where $c(G-S)$ denotes the number of components of $G-S$.
We present two tight lower bounds for $t(G)$ in terms of the Laplacian eigenvalues
and provide strong support for a conjecture for a better bound which, if true,
implies both bounds, and improves and generalizes known bounds by Alon, Brouwer, and the first author.
As applications, several new results on perfect matchings, factors and walks from Laplacian eigenvalues are obtained,
which leads to a conjecture about Hamiltonicity and Laplacian eigenvalues.
\end{abstract}



\section{Introduction}
Throughout this paper, $G=(V,E)$ is a simple graph of order $n$ with nonempty vertex set $V$ and nonempty edge set $E$.
The minimum degree of $G$ is denoted by $\delta$.
For a subset $S\subset V$, the subgraph of $G$ induced by $V\setminus S$ is denoted by $G-S$, and
$c(G-S)$ is the number of components of $G-S$.

The {\bf toughness} $t(G)$ of a graph $G$ is defined as $t(G)=\min\{\frac{|S|}{c(G-S)}\}$,
where the minimum is taken over all proper subsets $S\subset V$ such that $c(G-S)>1$.
By convention, a complete graph has infinite toughness.
For any real number $r\ge 0$, $G$ is {\bf $r$-tough} if $t(G)\ge r$.
The graph toughness was introduced by Chv\'atal \cite{Chva73} in 1973 and has been extensively studied since then.
It is closely related to graph structures, including cycles, matchings, factors, spanning trees, and others (see \cite{BaBS06}).

Toughness of regular graphs from eigenvalues of adjacency matrices has been well studied by, among others,
\cite{Alon95, Brou95, CiWo14, CiGu16, Gu21, Gu21b}. We use $\lambda_i :=\lambda_i(G)$ to denote the $i$th largest 
eigenvalue of the adjacency matrix of $G$, and let $\lambda =\max\{|\lambda_2|, |\lambda_n|\}$.
It was Alon~\cite{Alon95} who first showed that for any connected $d$-regular graph $G$, $t(G)>\frac{1}{3}(\frac{d^2}{d\lambda+\lambda^2}-1)$.
By this result, Alon showed that for every $t$ and $g$ there are $t$-tough graphs of girth strictly greater than $g$, which
strengthened a result of Bauer, Van den Heuvel and Schmeichel~\cite{BaVS95} who showed the same for $g=3$, and
thus disproved the pancyclic conjecture of Chv\'atal \cite{Chva73} in a strong sense.
Almost at the same time, Brouwer~\cite{Brou95} independently showed that
$t(G)>\frac{d}{\lambda}-2$ for any connected $d$-regular graph $G$. He also conjectured that $t(G) \ge\frac{d}{\lambda}-1$
in \cite{Brou95, Brou96}.
Recently Brouwer's conjecture has been confirmed by the first author~\cite{Gu21b}.
In this paper we consider arbitrary graphs and look for lower bounds on $t(G)$ in terms of the eigenvalues of the Laplacian matrix.

\medskip
The {\bf Laplacian matrix} $L$ (also called {\bf combinatorial Laplacian} or {\bf discrete Laplacian}) of a graph $G$,
is defined by $L=D-A$, where $D$ is the diagonal degree matrix and $A$ is the adjacency matrix of $G$.
Let $\mu_i:=\mu_i(G)$ denote the $i$th smallest eigenvalue of the Laplacian matrix of $G$.
Then $L$ is positive semi-definite and $\mu_1(G)=0$.
The second smallest Laplacian eigenvalue $\mu_2(G)$, is known as the {\bf algebraic connectivity} of $G$.
We have $\mu_2(G)=0$ if and only if $G$ is disconnected.
Moreover, if $\kappa$ is the vertex connectivity, then
\begin{equation}
\label{algcon}
\mu_2\leq\kappa\leq\delta.
\end{equation}
The complement $\G$ of $G$ has eigenvalues $\mu_1(\G)=0$ and $\mu_i(\G) = n-\mu_{n+2-i}(G)$ for $i=2,\ldots,n$.
Therefore $\mu_n(G)\leq n$ and $\mu_n(G)=n$ if and only if $\G$ is disconnected.
If $G$ is regular of degree $d$, then $L=dI-A$ and therefore $\mu_i=d-\lambda_i$ for $i=1,\ldots,n$.
For these and other properties of the Laplacian matrix and its eigenvalues we refer to \cite{BrHa12}.

\bigskip
The paper is organized as following. Our main results will be in the next section.
The tools and the proofs will be presented in Sections~\ref{tools} and \ref{proofs}.
In Section~\ref{appls}, we will show applications on perfect matchings, factors and walks.
Since toughness is closely related to cycle structures, we include a conjecture on Hamilton cycles
from Laplacian eigenvalues.

\section{Results}
Recently, the second author made the following conjecture~\cite{Haem20}.
\begin{conjecture}[Haemers]
\label{Haemconj}
\begin{equation}
\label{bd0}
t(G)\ge \frac{\mu_2}{\mu_n -\delta}.
\end{equation}
\end{conjecture}

For $d$-regular graphs, this conjecture implies that $t(G)\ge \frac{d-\lambda_2}{ -\lambda_n}$, which is stronger than Brouwer's conjecture.
The conjecture is supported by the following theorem and proposition.
The proofs will be given in Section~\ref{proofs}.

\begin{theorem}\label{mutou}
\begin{equation}
\label{bd1}
t(G)\ge \frac{\mu_n\mu_2}{n(\mu_n-\delta)},
\end{equation}
and
\begin{equation}
\label{bd2}
t(G)\ge \frac{\mu_2}{\mu_n - \mu_2}.
\end{equation}
\end{theorem}

\begin{proposition}\label{conjprop}
Let $S\subset V$ be such that $t(G)=|S|/c(G-S)$.
Then Conjecture~\ref{Haemconj} is true in each of the following cases.
\begin{itemize}
\item[\normalfont (i)]
The complement of $G$ is disconnected,
\item[\normalfont (ii)]
All connected components of $G-S$ are singletons (i.e. $n-|S|=c(G-S)$),
\item[\normalfont (iii)]
The union of some components of $G-S$ has order $\frac{1}{2}(n-|S|)$,
\item[\normalfont (iv)]
$c(G-S)=2$.
\end{itemize}
\end{proposition}

Since $\mu_n\leq n$ and $\delta\geq\mu_2$, the conjectured bound~(\ref{bd0}) implies (\ref{bd1}) and (\ref{bd2}).
Note that (\ref{bd0}) and (\ref{bd1}) coincide if $\mu_n=n$, that is, if the complement of $G$ is disconnected.
Therefore (i) of Propoition~\ref{conjprop} follows from (\ref{bd1}).

The three bounds coincide and are tight in case $G$ is the complete multipartite graph $K_{n_1,\ldots,n_m}$ ($1<m<n$).
Indeed, assume $n_1\geq\ldots\geq n_m$ then $t(G)=(n-n_1)/n_1$, $\mu_n=n$ and $\mu_2=\delta=n-n_1$.

The bounds (\ref{bd1}) and (\ref{bd2}) are incomparable.
For example, when $G$ is the Petersen graph, then $t(G)=4/3$ and $\mu_2=2$, $\mu_{10}=5$ and $\delta=3$,
and so (\ref{bd0}), (\ref{bd1}), and (\ref{bd2}) give $t(G)\geq 1$, $1/2$ and $2/3$, respectively.
The complement $\G$ satisfies $t(\G)=3$ and $\mu_2(\G)=5$, $\mu_{10}(\G)=8$ and $\delta=6$, and
the bounds (\ref{bd0}), (\ref{bd1}), and (\ref{bd2}) give $t(\G)\geq 5/2$, $2$ and $5/3$ respectively.

\section{Tools}\label{tools}

The following separation inequality from \cite{Haem95} provides a bridge between graph parameters and Laplacian eigenvalues.
It can also be found in \cite[Proposition 4.8.1]{BrHa12}.
\begin{theorem}[\cite{Haem95}]
Suppose that $X$ and $Y$ are two  disjoint subsets of $V$ such that there is no edge between $X$ and $Y$.
Then
\begin{equation}
\label{muxy}
\frac{|X| |Y|}{(n-|X|)(n-|Y|)} \le \left(\frac{\mu_n - \mu_2}{\mu_n + \mu_2}\right)^2.
\end{equation}
\end{theorem}

\par\medskip
By the above separation inequality, a simple yet useful proposition has been proved in \cite{GuLiu20}.
\begin{proposition}[\cite{GuLiu20}]
\label{22l}
Let $S\subset V$ such that $G-S$ is disconnected.
Let $X$ and $Y$ be disjoint vertex subsets of $V\setminus S$ such that $X\cup Y= V\setminus S$ with $|X|\le |Y|$. Then
\begin{equation}
\label{xupp}
|X| \le \frac{\mu_n - \mu_2}{2\mu_n}\cdot n,
\end{equation}
and
\begin{equation}
\label{ssiz}
|S| \ge \frac{2\mu_2}{\mu_n -\mu_2}\cdot |X|,
\end{equation}
with each equality holding  only when $|X|=|Y|$.
\end{proposition}

For review purpose, we include a proof, which was originally given in \cite{GuLiu20}.
\begin{proof}
By (\ref{muxy}), we have
\begin{equation}
\label{prodxy}
|X| |Y| \le \left(\frac{\mu_n - \mu_2}{\mu_n + \mu_2}\right)^2 (n-|X|)(n-|Y|),
\end{equation}
Let $\beta = \frac{\mu_n - \mu_2}{\mu_n + \mu_2}$. Then $0<\beta \le 1$, as $G$ contains at least one edge and $G-S$ is disconnected. It follows that
$$ |X|^2 \le |X|\cdot |Y| \le \beta^2 (n-|X|)(n-|Y|)\le \beta^2 (n-|X|)^2,$$
that is $$|X| \le \beta (n - |X|),$$ and hence
\begin{equation}
\label{betax}
|X|\le \frac{\beta n}{1+\beta}=\frac{\mu_n - \mu_2}{2\mu_n}\cdot n,
\end{equation}
with the equality holding  only when $|X|=|Y|$.

\par\medskip
Also, since $|Y| = n - |S| - |X|$, by (\ref{prodxy}), we have
$$ |X| (n - |S| - |X|) = |X|\cdot |Y| \le \beta^2 (n - |X|)(n - |Y|) = \beta^2 (n - |X|)(|S| + |X|),$$
implying that
\begin{equation}
\label{xnbo}
|X| n \le \left(\beta^2 (n - |X|) + |X| \right) \left(|S| + |X| \right)
= \left(\beta^2 n + (1-\beta^2)|X| \right) \left(|S| + |X| \right).
\end{equation}
By (\ref{betax}), we have
$$(1-\beta^2)|X|\le (1-\beta^2) \cdot \frac{\beta n}{1+\beta} = (\beta - \beta^2) n,$$
which, together with (\ref{xnbo}), implies that
$$|X|n \le \left(\beta^2 n + (\beta - \beta^2) n \right) \left(|S| + |X| \right) = \beta n \left(|S| + |X| \right),$$
and we have
$$ |X| \le \beta \left(|S| + |X| \right).$$
Hence,
\begin{equation*}
|S| \ge \frac{1 -\beta}{\beta} |X| = \frac{2\mu_2}{\mu_n - \mu_2}\cdot |X|.
\end{equation*}
Since (\ref{betax}) was utilized, the equality holds in (\ref{ssiz}) only when $|X|=|Y|$.
\end{proof}

Generalizing Hoffman's ratio bound,
the following bound for the cardinality of an independent set of an arbitrary graph
has been obtained independently in \cite{GoNe08} and in \cite{LuLT07} (see also~\cite{HaeHof}).
\begin{theorem}[Godsil and Newman~\cite{GoNe08}, Lu, Liu and Tian~\cite{LuLT07}]
\label{indupp}
If $U$ is an independent set of $G$, then $$|U|\le \frac{\mu_n - \delta}{\mu_n}\cdot n.$$
\end{theorem}

\section{Proofs}\label{proofs}
Throughout this section we take $S\subset V$ such that $t(G)=|S|/c(G-S)$, and put $c=c(G-S)$.
First we prove (\ref{bd1}) of Theorem~\ref{mutou}.
\begin{proof}
Clearly the vertex connectivity $\kappa$ of $G$ satisfies $\kappa\leq |S|$, so $|S|\geq \mu_2$ by (\ref{algcon}).
By taking a vertex in each component of $G-S$ we obtain an independent set of cardinality $c$.
Therefore Theorem~\ref{indupp} gives $c\leq n(\mu_n-\delta)/\mu_n$, and so
\[
t(G)\geq \frac{\mu_n\mu_2}{n(\mu_n-\delta)}.
\]
\\[-32pt]
\end{proof}
~\\
For convenience we continue with the proof of (ii) of Proposition~\ref{conjprop}.

\begin{proof}
If $n-|S| = c$, then $V\setminus S$ is an independent set of $G$.
By use of Theorem~\ref{indupp} and (\ref{algcon}) we have
\[
t(G) = \frac{|S|}{c} \ge \frac{n-c}{c} = \frac{n}{c}-1 \ge \frac{\mu_n}{\mu_n -\delta}-1 = \frac{\delta}{\mu_n -\delta}
\ge \frac{\mu_2}{\mu_n -\delta}.
\]
\\[-32pt]
\end{proof}
~\\
Next we prove (\ref{bd2}) of Theorem~\ref{mutou}.

\begin{proof}
Let $H_1, H_2, \ldots, H_c$ be the vertex sets of the components of $G-S$.
Without loss of generality, suppose that $|H_1|\le |H_2| \le \cdots \le |H_c|$.
Above we proved that (\ref{bd0}) and therefore (\ref{bd2}) holds if $H_1,\ldots,H_c$ are singletons.
Thus, we may assume that $n-|S| \ge c+1$.
We claim that $H_1, H_2, \ldots, H_c$ can be partitioned into two sets $X$ and $Y$ such that $|Y|\ge |X|\ge c/2$.
If $c$ is even, we can simply define $X =\bigcup_{1\le i\le \lfloor c/2\rfloor} H_i$ and $Y = (V\setminus S)\setminus X$.
Now we assume $c$ is odd.
If $|H_{(c-1)/2}|\ge 2$, then define $X =\bigcup_{1\le i\le (c-1)/2} H_i$ and $Y = (V\setminus S)\setminus X$ as needed.
The remaining case is $|H_1| = \cdots = |H_{(c-1)/2}|=1$.
We can define $X =\bigcup_{1\le i\le (c+1)/2} H_i$ and $Y = (V\setminus S)\setminus X$,
and we need to show that $|Y|\ge |X|\ge c/2$. If $|H_{(c+1)/2}| =1$,
then $|X| = \frac{c+1}{2}$ and $|Y| = n-|S|-|X| \ge \frac{c+1}{2}$, since $n-|S| \ge c+1$.
If $|H_{(c+1)/2}| \ge 2$, then $|X| = \frac{c-1}{2} + |H_{(c+1)/2}| \ge \frac{c-1}{2} +2 > \frac{c}{2}$
and $|Y| = \sum_{i > (c+1)/2} |H_i| \ge 2\cdot \frac{c-1}{2}=c-1\ge c/2$. Switch $X$ and $Y$ whenever needed to get $|Y|\ge |X|$.
\\[3pt]
It follows that $c\le 2|X|$.
Thus, by (\ref{ssiz}),
\begin{equation*}
t(G) = \frac{|S|}{c}\ge
\frac{2\mu_2}{\mu_n -\mu_2}\cdot \frac{|X|}{c} \ge \frac{\mu_2}{\mu_n -\mu_2}.
\end{equation*}
\\[-32pt]
\end{proof}
~\\
The last two proofs of this section deal with (iii) and (iv) of Proposition~\ref{conjprop}.

\begin{proof} (iii):
In this case $V\setminus S$ can be partitioned into two sets $X$ and $Y$ both having cardinality $\frac{1}{2}(n-|S|)$,
such that there are no edges between $X$ and $Y$.
We apply (\ref{ssiz}) and find $|S|\geq \frac{\mu_2}{\mu_n-\mu_2}(n-|S|)$, which implies
$|S|\geq  n\mu_2/\mu_n$.
As before, Theorem~\ref{indupp} gives $c\leq n(\mu_n-\delta)/\mu_n$ and hence
\[
t(G)= \frac{|S|}{c} \geq n\frac{\mu_2}{\mu_n} \cdot \frac{\mu_n}{n(\mu_n-\delta)}=\frac{\mu_2}{\mu_n-\delta}.
\]
\\[-32pt]
\end{proof}
\begin{proof} (iv):
It is known (see \cite{BrHa12}, Section 3.9) that $\mu_n\geq d_{max}+1$, when $d_{max}$ is the maximum degree of $G$.
If $G$ is not regular, then $d_{max}-\delta\geq 1$, and hence $\mu_n-\delta\geq 2$.
If $G$ is regular of degree $d=\delta$, then the adjacency matrix has smallest eigenvalue $\lambda_n=d-\mu_n$.
If $G$ is regular with $\lambda_n > -2$ then $G$ is the complete graph $K_n$ or an odd cycle $C_n$; see Theorem 2.5 of ~\cite{DoCv}.
We have $t(K_n)=\infty$ and $t(C_n)=2$ if $n\geq 4$.
If $n$ is odd, $\mu_2(C_n) = 2 - 2\cos(\pi/n)$ and $\mu_n(C_n) = 2 + 2\cos(2\pi/n)$, so (\ref{bd0})
gives $2\geq (1-\cos(\pi/n))/\cos(2\pi/n)$ which is clearly true for all odd $n\geq 5$.
Thus we can assume that $\lambda_n\leq -2$ and hence $\mu_n-d=\mu_n-\delta \geq 2$.
Thus we find
\[
t(G)=\frac{|S|}{c}\geq \frac{\mu_2}{2} \geq \frac{\mu_2}{\mu_n-\delta}.
\]
~\\[-32pt]
\end{proof}

\section{Applications}\label{appls}

It was shown in \cite{LiCh10} that if $\frac{\mu_2}{\mu_n} \ge \frac{2}{3}$, then $G$ is $2$-tough.
Now we can generalize it by rewriting (\ref{bd2}) of Theorem~\ref{mutou} as below.
\begin{theorem}
\label{toucor}
If $\displaystyle \frac{\mu_2}{\mu_n} \ge \frac{r}{r+1}$, then $G$ is $r$-tough.
\end{theorem}

Since many graph parameters and properties are related to toughness,
we have various applications, including but not limited to the results in this section.
We refer readers to the the survey paper~\cite{BaBS06} for more toughness related results.

The spectral conditions for matchings and $k$-factors of regular graphs have been attracting many researchers
\cite{BrHa05, Cioa05, KrSu06, CiGr07, CiGH09, OCi10, Lu10, Lu12, Gu14}, among others.
However, not as much has been discovered for general graphs that are not necessarily regular.
Brouwer and the second author~\cite{BrHa05} showed that if $n$ is even and  $2\mu_2 \ge \mu_n$, then $G$ has a perfect matching.
This result has been recently generalized to matching numbers in \cite{GuLiu20}.
We will have other generalizations by using graph toughness.

A graph $G$ is called {\bf elementary} if it contains a perfect matching and if the edges which occur in at least one perfect matching in $G$
induce a connected subgraph. A substantial study of elementary graphs has been given in \cite{LoPl86}.
It is proved in \cite{BBKMSS07} that every 1-tough graph with an even number of vertices is elementary.
Thus Theorem~\ref{toucor} implies the following result.

\begin{theorem}
If $n$ is even and $2\mu_2 \ge \mu_n$, then $G$ is elementary.
\end{theorem}

Let $G$ be an $n$-vertex graph with a perfect matching, and $m$ be a positive integer with $m<n/2 -1$.
Then $G$ is called {\bf $m$-extendable} if every matching of size $m$ extends to a perfect matching.
Plummer~\cite{Plum88} proved that if $t(G) > m$, then $G$ is $m$-extendable.

\begin{theorem}
Suppose $n$ is even, and let $m$ be a positive integer such that $m<n/2 -1$.
If $$\frac{\mu_2}{\mu_n} > \frac{m}{m+1},$$ then $G$ is $m$-extendable.
\end{theorem}

In \cite{EJKS85}, it is proved that every $k$-tough graph has a $k$-factor if $k|V(G)|$ is even and $|V(G)|\ge k+1$, confirming a conjecture
of Chv\'atal~\cite{Chva73}. The result was extended to non-regular factors by Katerinis~\cite{Kate90}.
Let $a\le b$ be positive integers. An {\bf $[a,b]$-factor} of a graph $G$ is a spanning subgraph $H$ such that $a\le d_H(v)\le b$ for each vertex.
Katerinis~\cite{Kate90} showed that for a graph $G$ on $n$ vertices such that $a<b$ or $bn$ is even, if $t(G)\ge a + \frac{a}{b} -1$,
then $G$ has an $[a,b]$-factor.
By using Theorem~\ref{toucor}, we have the following Laplacian eigenvalue condition for the existence of factors.

\begin{theorem}
Let $a\le b$ be positive integers such that $a<b$ or $bn$ is even.
If $$\frac{\mu_2}{\mu_n} \ge 1-\frac{b}{a(b+1)},$$ then $G$ has a $[a,b]$-factor.
\end{theorem}

The following corollary on $k$-factors is a generalization of the result on perfect matching by Brouwer and the second author~\cite{BrHa05}.
\begin{corollary}
\label{k-fac}
Let $k$ be a positive integer such that $n\ge k+1$ and $kn$ is even.
If $$\frac{\mu_2}{\mu_n} \ge \frac{k}{k+1},$$ then $G$ has a $k$-factor.
\end{corollary}

A graph $G$ is {\bf $(k,s)$-factor-critical} if $G-X$ has a $k$-factor for all $X\subseteq V(G)$ with $|X|=s$.
A $(1,1)$-factor-critical graph is usually referred to as a factor-critical graph.
It is proved in \cite{BBKMSS07} that every 1-tough graph with an odd number of vertices is factor-critical.
By Theorem~\ref{toucor}, we have the following result on factor-critical graphs from Laplacian eigenvalues,
which was originally obtained in \cite{GuLiu20}.
\begin{theorem}[\cite{GuLiu20}]
If $n$ is odd and $2\mu_2 \ge \mu_n$, then $G$ is $(1,1)$-factor-critical.
\end{theorem}

For $2\le s<n$, it was shown by Favaron~\cite{Fava96} that for a graph $G$ on $n$ vertices with $n+s$ even, if $t(G) > s/2$,
then $G$ is $(1,s)$-factor-critical.
By Theorem~\ref{toucor}, we have the following result.
\begin{theorem}
Suppose $2\le s<n$ and $n+s$ is even.
If $$\frac{\mu_2}{\mu_n} > \frac{s}{s+2},$$ then $G$ is $(1,s)$-factor-critical.
\end{theorem}
For $k=2, 3$ or even a general $k$, similar results on $(k,s)$-factor-critical graphs from toughness can be found in the survey~\cite{BaBS06}.
Thus Theorem~\ref{toucor} implies Laplacian eigenvalue conditions for $(k,s)$-factor-critical graphs with various values of $k$ and $s$,
which will be omitted here.

\medskip
Spectral conditions of the existence of a spanning tree with degree bounded above by a fixed number $k$ in a regular graph
have been obtained in \cite{CiWo14, CiGu16}.
When $k=2$, such a spanning tree is exactly a Hamilton path, and a sufficient condition has been given by, among others,
Butler and Chung~\cite{BuCh10}, who borrowed the idea from \cite{KrSu03} (both \cite{KrSu03, BuCh10} studied a stronger
structure, i.e., Hamilton cycle). The case of $k\ge 3$ for general graphs was asked in \cite{CiGu16} and has been solved in \cite{GuLiu20}.
A theorem of Win~\cite{Win89} implies that if $t(G)\ge \frac{1}{k-2}$ for $k\ge 3$, then $G$ has a spanning tree with maximum degree
at most $k$. Thus, Theorem~\ref{toucor} implies the following result of \cite{GuLiu20}.
\begin{theorem}[\cite{GuLiu20}]
\label{k-tre}
Let  $k\ge 3$ be an integer. 
If $$\frac{\mu_2}{\mu_n} \ge \frac{1}{k-1},$$ then $G$ has a spanning tree with maximum degree at most $k$.
\end{theorem}

Generalizing the idea of a Hamilton cycle, a {\bf $k$-walk} in a graph $G$ is a closed spanning walk of $G$ that visits every vertex of $G$
at most $k$ times. In particular, a Hamilton cycle is a 1-walk. Jackson and Wormald~\cite{JaWo90} observed that the existence of a spanning tree
with maximum degree at most $k$ actually implies the existence of a $k$-walk. Thus Theorem~\ref{k-tre} implies the existence of a
$k$-walk for $k\ge 3$. For $k=2$, Ellingham and Zha~\cite{ElZh00} showed that every 4-tough graph has a 2-walk.
By Theorem~\ref{toucor}, we have the following Laplacian eigenvalue condition for the existence of a $2$-walk.
\begin{theorem}
\label{2-wal}
If $\frac{\mu_2}{\mu_n} \ge \frac{4}{5}$, then $G$ has a $2$-walk.
\end{theorem}

\smallskip
\begin{remark}
In this section we presented applications of (\ref{bd2}) on perfect matchings, factors and walks.
The bound (\ref{bd1}) has similar applications and if Conjecture~\ref{Haemconj} is true,
then all results in this section can be improved in a similar manner.
\end{remark}

\begin{remark}
The results in this section as well as \cite{Haem95, BrHa05, YoLi12, GuLiu20} indicate that many graph properties are related to
the {\bf Laplacian eigenratio $\mu_2/\mu_n$}.
This eigenratio also has application aspects, and is highly related to the synchronization in Small-world Systems~\cite{BaPe02}.
It plays a similar role as the spectral gap, but may give a bit more information about the structure of the graph.
We feel that $\mu_2/\mu_n$ is interesting for future research.
\end{remark}

Similar to Theorems~\ref{k-tre} and \ref{2-wal}, the first author ever made the following conjecture for Hamilton cycles, but never
published elsewhere before.
\begin{conjecture}[Gu]
\label{conjham}
There exists a positive constant $C<1$ such that if $\mu_2/ \mu_n \ge C$ and $n\ge 3$ (or $n$ is sufficiently large),
then $G$ contains a Hamilton cycle.
\end{conjecture}

\begin{remark}
Notice that $K_{s, s+1}$ contains no Hamilton cycle, but $\mu_2 =s$ and $\mu_n = 2s+1$, which implies that
$\mu_2/ \mu_n$ can be arbitrarily close to $1/2$ if $s$ is sufficiently large. Thus the smallest possible $C$ is at least $1/2$.
\end{remark}

Chv\'atal~\cite{Chva73} conjectured that there exists a constant $t_0$ such that every $t_0$-tough graph contains a Hamilton cycle.
By (\ref{bd2}) of Theorem~\ref{mutou}, clearly Chv\'atal's conjecture implies Conjecture~\ref{conjham}.
Krivelevich and Sudakov~\cite{KrSu03} conjectured that for a $d$-regular graph $G$, there exists a constant $K$ such that
$d/\lambda > K$ and $n$ is large enough, then $G$ contains a Hamilton cycle.
It is not hard to see that Conjecture~\ref{conjham} implies the conjecture of Krivelevich and Sudakov.
All the three conjectures remain open.

\par\bigskip
\noindent
{\bf Acknowledgments}
\\The first author is partially supported by a grant from the Simons Foundation (522728, XG).


\end{document}